\documentclass[10pt,twocolumn,twoside]{IEEEtran}
\usepackage[utf8]{inputenc}
\newcommand{\comm}[1]{}
\usepackage{amssymb}
\usepackage{amsmath}
\usepackage{graphicx}
\usepackage{epstopdf}
\AppendGraphicsExtensions{.tif}

\usepackage[numbers]{natbib}\usepackage{epsfig}
\def\citet{\cite}
\def\iindex{}

\usepackage{times}\hfuzz=10pt 

 \usepackage{color}

\def\defi{\stackrel{{\scriptscriptstyle \Delta}}{=}}

\def\a{\alpha}

\def\w{\widehat}
\def\Ind{{\mathbb{I}}}

\def\R{{\bf R}}
\def\E{{\bf E}}
\def\P{{\bf P}}

\def\b{\beta}
\def\s{\delta}

\def\X{{\cal X}}
\def\t{\theta}
\def\oo{\bar}
\def\s{\sigma}

\def\G{\Gamma}

\newcommand{\be}{\begin{equation}}
\newcommand{\ee}{\end{equation}}
\newcommand{\bd}{\begin{displaymath}}
\newcommand{\ed}{\end{displaymath}}
\newcommand{\ba}{\begin{array}{ll}}
\newcommand{\ea}{\end{array}}
\newcommand{\baa}{\begin{eqnarray}}
\newcommand{\eaa}{\end{eqnarray}}
\newcommand{\baaa}{\begin{eqnarray*}}
\newcommand{\eaaa}{\end{eqnarray*}}

\def\oo{\bar}

\def\a{\alpha}

\def\ee{\epsilon}

\def\EE{{\mathbb{E}}}
\def\VVar{{\mathbb{V}ar}}

\title{On approximation of the distribution for Pearson statistic}
\author{
Nikolai Dokuchaev}

\begin{document}
 \vspace{-0.5cm}
\def\brea{}
\def\breakk{}
\def\break{}
\def\break{\nonumber\\ }\def\breakk{\nonumber\\&&}\def\brea{\nonumber\\ }
\maketitle 
 \let\thefootnote\relax\footnote{The author is with the School of Electrical Engineering, Computing and Mathematical Sciences, Curtin University, GPO Box U1987, Perth,
Western Australia, 6845 \comm{and  ITMO University, St. Petersburg, Russian Federation, 197101.}}
\begin{abstract} The paper considers the classical Goodness of Fit test.
 It suggests to use the Gamma distribution for the approximation of the distribution of the Pearson statistics with unknown parameters estimated from raw data. The parameters of these Gamma distribution can be estimated from  the first two moments of the  statistic after  averaging over a distribution of the unknown  parameter over its range. This allows  to
simplify calculation of the quantiles for the Pearson statistic, as is shown in some simulation experiments  with medium and small sample sizes.
\end{abstract}

\section{Introduction}
The classical statistical  Goodness of Fit  test  addresses the problem
of estimation the parameters of a parametric family of distributions from observed data with unknown $d$-dimensional parameter that has to be fitted from the data.
The Pearson statistics  is commonly used to estimate the error; see, e.g., the literature review in \cite{Ba,GN,P}. Let $n$ be the number in intervals where the observations are counted
in the Pearson statistic.
The limit distribution of this statistics for infinitely increasing sample size is known, given some mild conditions; see, e.g. \cite{Ba,B,GN,P}.  The quantiles
for its limit distribution  are often used as the critical values for the test.
If the parameter is fitted from the raw (ungrouped)  data using a consistent estimator, then the limit distribution is different; see, e.g. \cite{Ba}, p.24.
The actual distribution of statistic for the finite samples is
a discreet distribution and depends on the choice of the counting intervals and other  parameters of the experiment.

 This sort paper  suggests to use the Gamma distribution for a simplified  approximation of the distribution of the Pearson statistic with small and medium sample sizes.
 The parameters of these Gamma distribution can be estimated from  the first two moments of the sample distribution
 of simulated Pearson statistic with  parameter values randomized  over a domain for the unknown true parameter.
Some  computer experiments with medium and small sample sizes shows that this helps to reduce the bias
for the calculation of the quantiles for the Pearson statistic.

\section{Problem setting}
Let $F(\cdot|\t)|_{\t\in D}$ be a give family of distributions, where $\t\in D$ is a $d$-dimensional parameter, and where $D\subset \R^d$ is a domain.

Assume that we are testing a   hypothesis  about a population distribution
for a given independent and identically distributed sample $X=(X_1,...,X_N)$
from the distribution $F(\cdot|\t_0)$, where $\t_0\in D$. Let $\w\t$ be the estimate of $\t$ obtained using a consistent estimator $\w\t=T(X)$, where
$T:\R^N\to \R^d$ is  a mapping. We assume that $\w\t$ is observable.

For a given integer $n>0$, consider a system of mutually disjoint intervals $\{I_i\}_{k=1}^n$ such that
$\cup_i I_i=\R$ (two of these intervals are semi-infinite).
Let $C_i$ be observed sample counts in the intervals
$I_i$, calculated from a sample $x_1,...,x_N$. In particular,  we have  that
\baaa\sum_{i=1}^nC_i=N.\eaaa

The values $C_i$ are supposed to be observable.

Consider two hypotheses:

$H_0$: The sample come from the distribution $F(\cdot|\w\t)$.

$H_A$: The sample does not come from this distribution.

A hypothesis has to be accepted or rejected based on observed  $\w\t$ and observed
counts $\{C_i\}$, given that the family of the distributions $F(\cdot|\t)|_{\t\in D}$, and the domain $D$ of possible values of $\t$ are known.

Let $ Q_i\defi \E\{C_i|\t=\w\t\}$.

Let us consider Pearson's statistic  \baaa \X^2\defi\sum_{i=1}^n \frac{(C_i-Q_i)^2}{Q_i}. \eaaa
If
the computed value of $\X^2$ is large, then we reject hypothesis $H_0$.
In this case,  the observed and expected values are not close and
the model is a poor fit to the data.

  If the parameter is fitted from the grouped data using a consistent estimator based on counting in the intervals, then the limit distribution
is a known  $\chi^2_{n-d-1}$-distribution, given some mild conditions; see, e.g. \cite{Ba,B,GN,P}.
If the parameter is fitted from the raw (ungrouped)  data using a consistent estimator, then the limit distribution
is \baa
\chi^2_{n-d-1}+\sum_{k=n-d}^{n-1}\nu_k Z_k,
\eaa where $Z_k$ are independent  standard normal variables
and $\nu_k\in[0,1]$ (Chernoff and Lehmann \cite{CL}; see also \cite{Ba}, p.24).
However, the values $\{\nu_k\}$ depend on the intervals, on the population distribution, and
on the estimator.

The distribution of $\X^2$ is discreet and depends on the choice of
$
(F(\cdot|\cdot),D,\t_0,\{I_i\}_{i=1}^n,T(\cdot))$.
The standard approach for the approximation of the distribution of $\X^2$ for large $N$
is its approximation  $\chi^2_{n-1-d}$ distribution, i.e,,  by the $\chi^2$-distribution with $n-1-d$ degrees of freedom (see, e.g., \cite{B}).  In the literature, $X^2$ is called {\em chi square
statistic} or {\em Pearson's statistic}.  This  $\chi^2_{n-1-d}$ distribution is independent on the choice of the set $\{I_i\}$. However, the actual distribution of $\X^2$ is not easy to describe; it
depends on   $(F(\cdot|\cdot),D,\t_0,\{I_i\}\{I_i\}_{i=1}^n,T(\cdot))$. Therefore it is not easy to calculate quantiles used for the hypothesis testing.  On the other hand, some numerical examples given below show that use of  quantiles for the $\chi^2_{n-1-d}$ distribution as a substitution for quantiles  of $X^2$ could lead to significant bias for the critical values.

\iindex{The distribution of $\X^2$ is discreet and depends on the choice of $\{n,\{I_i\}_{i=1}^n,F(\cdot|\t)\}$.
It is known  that, under some mild assumptions,  \baa
&&\hbox{the distribution of}\quad  \X^2\quad \hbox{converges to} \quad   \chi^2_{n-1-d}\breakk \quad \hbox{as}\quad  N\to +\infty.\label{conv}\eaa
(See, e.g., \cite{B}). The limit distribution here  is the
$\chi^2$-distribution with $n-1-d$ degrees of freedom. In the literature, $X^2$ is called {\em chi square
statistic} or {\em Pearson's statistic}.  This limit distribution is independent on  $
(F(\cdot|\cdot),D,\t_0,\{I_i\},T(\cdot))$. However, the actual distribution of $\X^2$ is not easy to describe for a given finite $n$. Therefore it is not easy to calculate quantiles for the hypothesis testing.  On the other hand, some numerical examples given below show that use of  quantiles for the $\chi^2_{n-1-d}$ distribution as a substitution for quantiles  of $\X^2$ could lead to significant bias for the critical values.}

There are several known approaches to deal with this bias; see, e.g. \cite{GN}.
We suggest one more approach that seems to provide a reasonably close approximation
for the distribution of the tests statistics with medium and small sample sizes.
\section{Approximation by Gamma distribution}
In some numerical experiments, we have found that  the Gamma distribution $\G(\a,\lambda)$ with the density
$\Ind_{\{x>0\}}\frac{\lambda^\a}{\G(\a)}x^{\a-1}e^{-\lambda x}$  can be effectively  used as a close approximation of the distribution of $\X^2$.

For the  distribution  $\G(\a,\lambda)$, the expectation is $\a/\lambda$,
and the variance is $\a^2/\lambda$.

Technically, the distribution of $\X^2$ and as well as $(\a,\lambda)$ depend on the choice of
$(F(\cdot|\cdot),D,\t_0,\{I_i\}_{i=1}^n,T(\cdot))$.  However, we  need an approximation that  does not use $\t_0$.   Therefore, we suggest to estimate these parameters via matching them
 with the first two moments of a sample of random values  $\X^2$ simulated  under the compound  distribution $F(\cdot|\Theta)$ with a random $\Theta$ given some preselected probability distribution for $\Theta$. This removes dependence of $(\a,\lambda)$ on the true parameter $\t$. For example, one can use a non-informative uniform distribution over a bounded domain $D$ containing the true parameter  $\t$
 \par
Let $\EE$ and $\VVar$ be the sample mean and the sample variance, respectively,
over the Monte-Carlo trials for simulation of  $(\Theta,X)$ generating the implied statistic $\X^2$.

  \subsubsection*{The procedure for fitting $(\a,\lambda)$}
 \begin{enumerate}
 \item  Run $M$ Monte-Carlo simulations of $\Theta$. For each simulated $\Theta$, simulate an i.i.d. sample  $X=(X_1,...,X_N)$ with the terms
 distributed under  $F(\cdot |\Theta)$.
 \item Calculate $\X^2$ for each simulation of $(\Theta,X)$.
 \item Calculate $a=\EE \X^2$ and $v=\VVar X^2$.
 \item Find $\a$ and $\lambda $ such that
 \baaa
 \frac{\a}{\lambda}=\EE \X^2,\quad \frac{\a}{\lambda^2}=\VVar \X^2,
 \eaaa
 i.e.
\baaa
\a= \frac{(\EE \X^2)^2}{\VVar \X^2},\quad  \lambda= \frac{\EE \X^2}{\VVar \X^2}.
 \eaaa
 \item Use quantiles for $\G(\a,\lambda)$ as approximations for quantiles for $\X^2$.
  \end{enumerate}
It seems that this approach allows to achieve a significant reduction of the bias for quantiles  for the sample sizes.\comm{ that are insufficient to ensure good approximation  by the limit distribution; these sizes include $N\le 2000$, and possibly even larger $N$.}

 \section{Numerical examples}
\label{SecN}
Let illustrate the difference between the limit distribution and actual distribution of $\X^2$ using the following numerical example.

This would correspond to the setting with $d=1$  \comm{ in (\ref{conv})} and $n-d-1=1$. .

We run Monte-Carlo experiments with  the sample size $M=10^6$ for $\X^2$.
We run these experiments for four cases with different sets  of parameters. These cases are  listed below.

{\bf Case A:}

For this case, we simulated $\X^2$ for the sample $X$
from exponential distribution $Exp(\t_0)$, i.e. with the density
$\Ind_{\{x>0\}}\t^{-1}e^{-\t x}$. This corresponds to the case  of non-random $\Theta=\t_0$. We have  used $\t_0=1$, and we have used the estimate $\w\t=T(X)=1/\oo X$, where $\oo X\defi \frac{1}{N}\sum_{i=1}^N X_i$.
This is a maximum likelihood estimate as well as the estimate implied by the method of moments. It is known that this estimate is consistent.

{\em Case A(i)}:
The sample size for the underlying process $X$  is
$N=10$, the number of intervals  is $n=3$, and
 $I_1=(-\infty,a_1]$, $I_2=(a_1,a_2)$, $I_3=[a_2,\infty)$. The numbers  $a_1<a_2$
are such that
$\P(X_k\in I_k|\t_0)=1/3$. This choice corresponds to the most basic  case where of equal probabilities for the intervals.
For this case, we found in the experiments with $10^6$ Monte-Carlo trials that $\EE \X^2= 1.35800$ and $\VVar  \X^2=2.0845822$.

As can be seen, it is quite far for the  expectation
 and the variance for the $\chi^2_{n-d-1}=\chi_1^2$ distribution; these parameters are $n-d-1=1$ and $2(n-d-1)=2$
respectively.

{\em Case A(ii)}: The sample size for the underlying process $X$  is
$N=1000$; the remaining parameters are the same as for Case A(i).
For Case A(ii), we found in the experiments with $10^6$ Monte-Carlo trials that $\EE \X^2=1.350675$ and $\VVar  \X^2=2.245898$.

Table \ref{tabA} shows sample quantiles for $\X^2$ for Cases A(i)-(ii). This example shows that use of  quantiles for the limit distribution as a substitution for quantiles  of $X^2$ for finite samples could lead to a
bias. Approximation by Gamma function helps to reduce the bias,
as is shown in examples described below.

\begin{table}[ht]
 \begin{center}  (i) Quantiles for Case A(i); the sample size for $X$ is $N=10$\\
 \begin{tabular}{|c|c|c|c|c|}\hline
 {\em Quantiles}   &{\em 0.75} & \em {0.9} & {\em 0.95} & {\em 0.99}
 \\
\hline For  $\X^2$  &1.801390 &3.052967 &4.146487 & 6.279877                    \\
\hline
 \end{tabular}
\end{center}
 \begin{center} (ii) Quantiles for Case A(ii); the sample size for $X$ is $N=1000$\\
 \begin{tabular}{|c|c|c|c|c|}\hline
 {\em Quantiles}   &{\em 0.75} & \em {0.9} & {\em 0.95} & {\em 0.99}
 \\
\hline For  $\X^2$  &  1.810692 & 3.195782 & 4.312670 & 7.092106 \\
\hline
 \end{tabular}
\end{center}
\caption{Quantiles for Cases A(i)-A(ii).}
\label{tabA}\end{table}

We have also considered  cases where  the parameters $(\a,\lambda)$  have been fitted to the sample $X^2$ simulated  according to the procedure describes above.

{\bf Case B:} For this case, we  consider the family the exponential distribution $Exp(\t)$ with the density
$\Ind_{\{x>0\}}\t^{-1}e^{-\t x}$.  We  assumed that   $\t\in D=[0.5,1.5]$. For the step (i) of this  procedure, we have used  $\Theta$ uniformly  distributed
on the domain $D=[0.2,2]$. Further,  the sample size  for the underlying process $X$ for this case is
$N=20$, the number of intervals  is $n=3$, and
 $I_1=(-\infty,0.5]$, $I_2=(0.5,1.5)$, and $I_3=[1.5,\infty)$.

 For this case, we have $\EE \X^2= 1.323495$, $\VVar \X^2=2.142576$, and the corresponding
 parameters for $\G(\a,\lambda)$ are
\baaa
 \a= 0.8175386,\quad \lambda=0.617712.
 \eaaa
 Table \ref{tab1}(i) shows quantiles for $\X^2$, for the fitted distribution $\G(\a,\lambda)$, and
 for the limit distribution  $\chi^2_{n-d-1}$.

{\bf Case C:} For this case, we consider $I_1=(-\infty,1]$, $I_2=(1,2)$, and $I_3=[2,\infty)$.
All other parameters are the same as for Case B.

 For this case, we have $\EE \X^2=1.294582$, $\VVar \X^2=2.183124$, and the corresponding
 parameters for $\G(\a,\lambda)$ are \baaa
 \a=0.7676803,\quad \lambda=0.5929949
 \eaaa
 Table \ref{tab1}(ii) shows quantiles for $\X^2$, for the fitted distribution $\G(\a,\lambda)$, and
 for the limit distribution  $\chi^2_{n-d-1}$.

{\bf Case D:}  For this case, we consider
$N=1000$.
All other parameters are the same as for Case C.

 For this case, we have $\EE \X^2=1.281905$, $\VVar \X^2=2.191206$, and the corresponding
 parameters for $\G(\a,\lambda)$ are \baaa
 \a= 0.7499429,\quad \lambda=  0.5850224.
 \eaaa
 Table \ref{tab1}(iii)  shows quantiles for $\X^2$, for the fitted distribution $\G(\a,\lambda)$, and
 for the limit distribution  $\chi^2_{n-d-1}$.

{\bf Case E:}  For this case,we  consider the family the normal distributions $N(\mu,\s^2)$ with
$\mu\in[-0.5,0.5]$ and $\s\in [1,2]$ and $\theta=(\mu,\s)$.
The random parameter $\Theta$ as a random vector with independent components
 distributed uniformly on $[-0.5,0.5]$ and $[1,2]$ respectively.
 We used $(\w\mu,\w\s)=T(X)$ such that $\w\mu$ is the sample mean of $X$ and
 $\w\s^2$ is the sample variance of $X$.
 The number of intervals  is $n=4$, and the intervals are
 $I_1=(-\infty,-1]$, $I_2=(-1,0]$,  $I_3=(0,1]$, and $I_4=(1,\infty)$.

 For this case, we have $\EE \X^2=1.772562$, $\VVar \X^2=2.852873$, and the corresponding
 parameters for $\G(\a,\lambda)$ are \baaa
 \a= 1.101338,\quad \lambda= 0.621325 .
 \eaaa
 Table \ref{tab1} (iv)   shows quantiles for $\X^2$ and for the fitted distribution $\G(\a,\lambda)$.

{\bf Case F:}  For this case,we  consider the family the normal distributions $N(\mu,\s^2)$ with
$\mu\in[-1,1]$ and $\s\in [0.5,4]$ and $\theta=(\mu,\s)$.
The random parameter $\Theta$ as a random vector with independent components
 distributed uniformly on $[-1,1]$ and $[0.5,4]$ respectively.
 The intervals and the estimatres are the same as in Case E.
 $I_1=(-\infty,-1]$, $I_2=(-1,0]$,  $I_3=(0,1]$, and $I_4=(1,\infty)$.

 For this case, we have $\EE \X^2= 1.922529$, $\VVar \X^2=3.251841$, and the corresponding
 parameters for $\G(\a,\lambda)$ are \baaa
 \a=1.136623,\quad \lambda= 0.5912125
.
 \eaaa
 Table \ref{tab1} (v)   shows quantiles for $\X^2$ and for the fitted distribution $\G(\a,\lambda)$.
[

\begin{figure}[ht]
\centerline{\psfig{figure=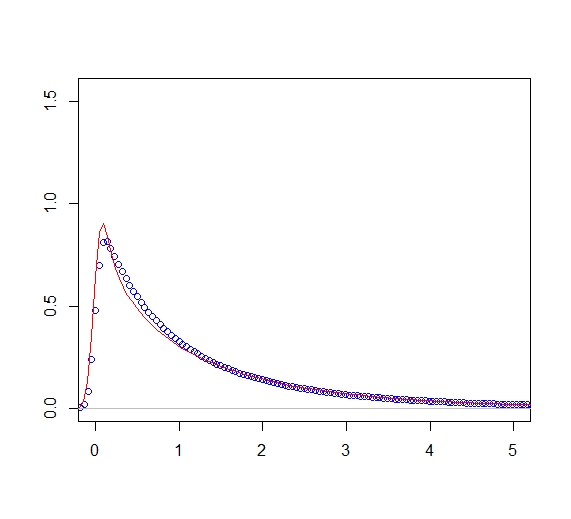,width=9cm,height=5.0cm}} \centerline{\psfig{figure=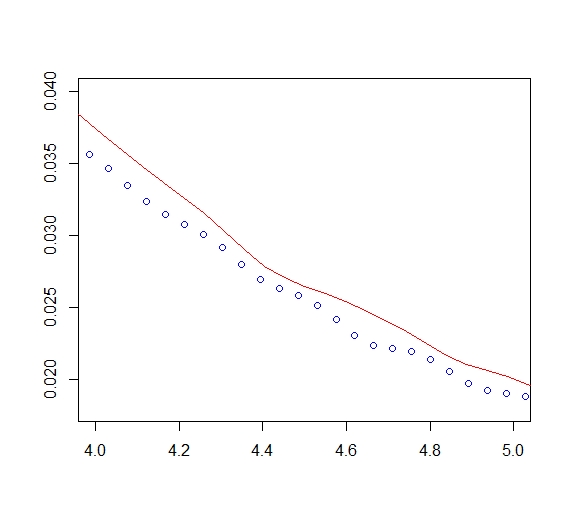,width=9cm,height=4.0cm}}
\caption[]{ Smoothed histograms  for $\X^2$ (circles) and $\G(\a,\lambda)$ (line),
 recovered from $10^6$-size simulated sample  for the Case D, in two different
 magnifications.}
\label{fig1}
\end{figure}

\begin{figure}[ht]
\centerline{\psfig{figure=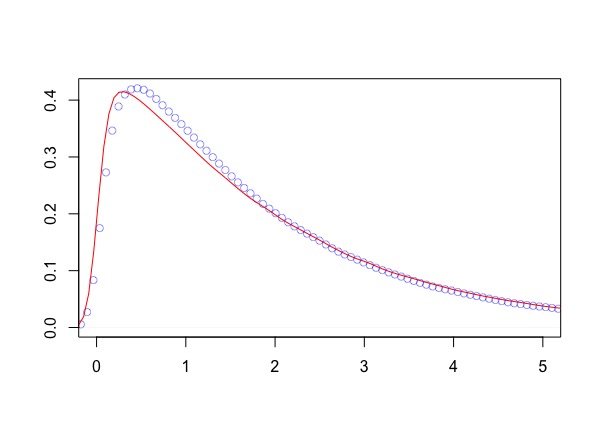,width=9cm,height=5.0cm}} \centerline{\psfig{figure=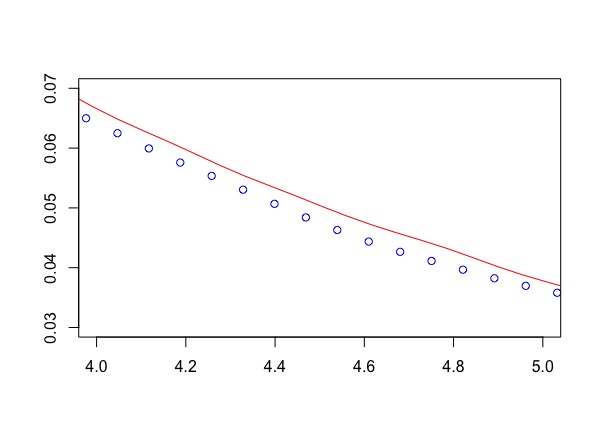,width=9cm,height=4.0cm}}
\caption[]{ Smoothed histograms  for $\X^2$ (circles),= and $\G(\a,\lambda)$ (line),
 recovered from $10^6$-size simulated sample  for the Case F, in two different magnifications.}
\label{fig2}
\end{figure}

Figures \ref{fig1}-\ref{fig2}  show smoothed histograms   for $\X^2$ and $\G(\a,\lambda)$
for Cased D,E,and F, respectively, constructed  from the histograms for  Monte-Carlo samples of the size  $M=10^6$  using the standard command {\em densities} in {\em R} programming language.    These figures demonstrate quite close approximation.

We have used {\em R} programming language for calculations; calculation of $(\a,\lambda)$ for $N=20$ and $M=10^6$ takes less than a minute for a standard desktop computer. For $N=1000$ and $M=10^6$, it takes about 10 minutes.

\section{Conclusion}
The paper suggest to approximate the distribution of the Pearson statistic
 by the Gamma distributions with parameters fitted to simulated statistics a given configuration of cells  where the sample occurrences are being counted. Feasibility of this approach is demonstrate with some numerical experiments.   So far, the range of the parameters for these experiment was quite limited. It would be interesting to extend these experiments on more general choices of the parameters, especially  $n$ and $N$. We leave this for the future research.

\newpage
\begin{table}[ht]
 \begin{center}  (i) Quantiles for Case B; the sample size for $X$ is $N=20$\\
 \begin{tabular}{|c|c|c|c|c|}\hline
 {\em Quantiles}   &{\em 0.75} & \em {0.9} & {\em 0.95} & {\em 0.99}
 \\
\hline For  $\X^2$  & 1.787257 & 3.111514 & 4.296282 & 6.762272\\
   \hline  For fitted $\G(\a,\b)$ & 1.831157 & 3.204561 & 4.262158 & 6.760412  \\
\hline
 \end{tabular}
\end{center}
 \begin{center} (ii) Quantiles for Case C; the sample size for $X$ is $N=20$\\
 \begin{tabular}{|c|c|c|c|c|}\hline
 {\em Quantiles}   &{\em 0.75} & \em {0.9} & {\em 0.95} & {\em 0.99}
 \\
\hline For  $\X^2$  &1.690018 & 2.993245 & 4.200237 & 6.907514 \\
   \hline  For fitted $\G(\a,\b)$  & 1.781864 & 3.178560 & 4.255779 &6.811008  \\
  \hline
 \end{tabular}
\end{center}
 \begin{center} (iii) Quantiles for Case D; the sample size for $X$ is $N=1000$\\
 \begin{tabular}{|c|c|c|c|c|}\hline
 {\em Quantiles}   &{\em 0.75} & \em {0.9} & {\em 0.95} & {\em 0.99}
 \\
\hline For   $\X^2$  &  1.707008 &3.094988 & 4.231155 & 7.007469  \\                        
   \hline  For fitted $\G(\a,\b)$  &  1.765294 & 3.161543 & 4.250052 & 6.850396 \   \\
\hline
 \end{tabular}
 \end{center}
\begin{center} (iv) Quantiles for Case E; the sample size for $X$ is $N=1000$\\
 \begin{tabular}{|c|c|c|c|c|}\hline
 {\em Quantiles}   &{\em 0.75} & \em {0.9} & {\em 0.95} & {\em 0.99}
 \\
\hline For   $\X^2$  & 2.397755 & 3.939819 & 5.124125 & 7.927041  \\
   \hline  For fitted $\G(\a,\b)$ &  2.453660 & 3.982921& 5.121967 &7.796123    \\
\hline
 \end{tabular}
 \end{center}
 \begin{center} (v) Quantiles for Case F; the sample size for $X$ is $N=1000$\\
 \begin{tabular}{|c|c|c|c|c|}\hline
 {\em Quantiles}   &{\em 0.75} & \em {0.9} & {\em 0.95} & {\em 0.99}
 \\
\hline For   $\X^2$  &  2.602056 &4.228039 &5.479738 &8.478611   \\
   \hline  For fitted $\G(\a,\b)$ & 2.656888 & 4.284709 & 5.499895 & 8.278232     \\
\hline
 \end{tabular}
 \end{center}
 \caption{Quantiles for Cases B,C,D.}
\label{tab1}\end{table}
 \newpage
\end{document}